\documentclass[12pt,twoside]{myarticle}
\usepackage[all]{xy}

\usepackage{theorem}
\usepackage{epsfig}
\usepackage{amsfonts}
\usepackage{amssymb}

\theoremstyle{change}{\theorembodyfont{\slshape}
\theoremheaderfont{\scshape}
\newtheorem{theorem}{Theorem.}[section]
\newtheorem{proposition}[theorem]{Proposition.}
\newtheorem{lemma}[theorem]{Lemma.}
\newtheorem{remark}[theorem]{Remark.}
\newtheorem{corollary}[theorem]{Corollary.}

}

\theoremstyle{change}{
\theorembodyfont{\rmfamily}
\newtheorem{example}[theorem]{Example.}

\newtheorem{definition}[theorem]{Definition.}
}

\def\proof{\noindent{\sc Proof.}\enspace}
\def\endproof{ \quad $\kasten$\bigskip}

\def\S{\mathop{\cal{S}{}}}

\def\F{\mathop{\cal{F}{}}}
\def\C{\mathop{\cal{C}{}}}

\def\epsilon{\varepsilon}

\def\t#1{\widetilde{#1}}
\def\b#1{\overline{#1}}

\def\CC{{\mathbb C}}
\def\ZZ{{\mathbb Z}}
\def\RR{{\mathbb R}}
\def\NN{{\mathbb N}}

\def\PP{{\mathbb P}}
\def\KK{{\mathbb K}}

\def\good{\mathop{\rm good}\nolimits}

\def\good{/ \! /}

\def\mal{\mathbin{\! \cdot \!}}

\def\cone{\mathop{\hbox{\rm cone}}}
\def\lin{\mathop{\hbox{\rm lin}}}
\def\id{\mathop{\rm id}\nolimits}
\def\pr{\mathop{\rm pr}\nolimits}

\def\osubset{\subset \kern-10pt {\rm o}\ } 

\def\kasten{\mathord{\vbox{\hrule
                     \hbox{\vrule
                     \hskip5pt
                     \vrule height5pt
                     \vrule}
                     \hrule}}}
\def\text#1{\hbox{\rm #1}}

\def\ltexindent#1{\hbox to \hangindent{#1\hss}\ignorespaces}

\begin{document}

\title{Weakly Proper Toric Quotients%
\footnote{\small {\it Mathematics Subject Classification.}\enspace 
14M25, 14L30, 14D25.}}

\author{Annette A'Campo--Neuen
  \\
{e-mail: annette@math-lab.unibas.ch}
}

\maketitle

\pagestyle{myheadings}
\markboth{\hfill A.~A'Campo--Neuen\hfill}{\hfill Weakly proper toric quotients
\hfill}

\begin{abstract}

\noindent We consider subtorus actions on complex
toric varieties. A natural candidate
for a categorical quotient of such an action is the so--called {\it toric
quotient\/}, a universal object constructed in the toric category.
We prove that if the toric quotient is {\it weakly proper\/} 
and if in addition the quotient variety is of expected dimension
then the toric quotient is a categorical quotient in the category of algebraic
varieties. 
For example, weak properness always holds  
for the toric quotient of a subtorus action on a toric variety whose
fan has a convex support. 
\end{abstract}

\section*{Introduction}

\noindent
In \cite{GIT} D. Mumford introduced the notion of a categorical quotient
for the  action of an algebraic group $G$ on an algebraic variety $X$. 
By definition this is a $G$--invariant morphism $p\colon X\to Y$ such that
every $G$--invariant morphism from $X$ to some algebraic variety
factors uniquely through $p$.
In general, such a categorical quotient need not exist. In this article
we will consider subtorus actions on complex toric varieties.

In this setting, a natural candidate for a categorical quotient has
been constructed in the category of toric varieties, namely the
so--called {\it toric quotient} (see \cite{Torische Quotienten}).
The toric quotient is universal for  toric morphisms from the given toric
variety that are constant on the orbits of the subtorus action.
Clearly, a necessary condition for the toric quotient to be categorical
is surjectivity. But there are examples of toric quotients that are
not surjective and hence not categorical. In fact, in \cite{Beispielsammlung},
Section~5,
an example of a subtorus action on a toric variety is given that
does not admit a categorical quotient, not even in the category
of algebraic or analytic spaces.

On the other hand, if the codimension of $H$ in the big torus $T$
of the toric variety is at most $2$, then the toric quotient is always
categorical (Corollary 4.3 in \cite{Beispielsammlung}). 
An important tool for the proof
was to observe that the toric quotient in that case is {\it weakly
proper\/}, i.e.~it satisfies a certain weak
lifting property for holomorphic germs of curves (the precise definition is given in Section~1). 

A morphism of complex algebraic varieties
is weakly proper if and only if it is universally submersive.
Another more geometric characterization  is given
by the following notion:
We say that a morphism $p\colon X\to Y$ of algebraic varieties
satisfies the {\it curve covering 
property\/} if for every curve $Y'\subset Y$ 
there is a (not necessarily irreducible) 
curve $X'\subset X$ with $p(X')=Y'$,
such that the morphism from $X'$ to $Y'$ defined by $p$ has finite fibers.
If in addition the curve covering property is preserved under
any base change, then the morphism $p$ has
the {\it universal curve covering property\/}. If $p$ is a surjective
morphism over the field of complex numbers then this property
is equivalent to weak properness (see Section~1).

For a toric morphism, weak properness has a very
simple characterization in terms of fans, namely the property
holds if and only if
the associated lattice homomorphism induces a surjective map on
the supports of the corresponding fans (see Section~2).

In this article we will show that weak properness
is a sufficient condition for a toric quotient of expected dimension
to be categorical. More precisely, our main result is the following
(see Corollary \ref{Hauptergebnis}):

\bigskip

\noindent{\sc Theorem.\enspace} 
{\sl For a toric variety $X$ and a subtorus $H$ of the big torus $T$ of $X$, 
let $p\colon X\to Y$ denote
the toric quotient for the action of $H$ on $X$. 
If $p$ is weakly proper and $\dim Y=\dim T/H$,
then the toric quotient is categorical.}

\bigskip

For example, weak properness is automatically satisfied
if the fan associated to $X$ has a convex support, or equivalently,
if there is a proper toric morphism from $X$ onto an affine toric variety. In fact,
we can even show that in this case the dimension condition can be 
omitted (see Corollary \ref{convex case}):

\bigskip

\noindent{\sc Theorem.\enspace} 
{\sl If the toric variety $X$ corresponds to a fan with convex support,
then for any subtorus action on $X$ the toric quotient is a
quotient in the category of algebraic varieties.}

\bigskip

In order to prove our results we proceed by induction on the number of steps
occuring in the construction of the toric quotient. The intermediate steps
of this construction can be viewed as successive approximations of the quotient by
non-separated prevarieties (see Section~6). Therefore it simplifies the arguments
to  work in the more general context of toric prevarieties. We prove a result
stated in this context (see Theorem~\ref{weakly proper toric quotient}) and obtain
the above theorems as corollaries.

Toric prevarieties are the non--separated analogues of toric varieties, 
i.e.~complex algebraic prevarieties
with an effective action of a torus having a dense orbit.
In analogy to the separated case, there is a convex--geome\-tri\-cal
description of toric prevarieties in terms of so--called
{\it systems of fans} (see \cite{Prevarieties}).
In Section~2 we briefly recall the basic facts on toric prevarieties
and systems of fans. 

In Section~3 we consider morphisms from
toric prevarieties to algebraic varieties that are not necessarily
toric. Every such morphism defines
an equivalence relation on the cones occuring in the 
system of fans corresponding to the toric prevariety, and 
the supports of the equivalence classes form a partition
of the support of the system of fans. The properties of this partition
are analyzed in
Section~4 and~5. These sections contain the
convex--geometrical lemmata that are needed for the proof of the main theorem
that is carried out in Section~6.

As an application of the result, in Section~7
we give an example of a categorical quotient 
$p\colon X\to Y$ of
a $4$--dimensional toric variety by some $\CC^*$--action where
$\dim Y=3$ that is
not uniform in the sense of \cite{GIT}, i.e. such that for some open
subset $U\subset Y$, the restriction
$p\colon p^{-1}(U)\to U$  is not
the categorical quotient for the induced $\CC^*$--action.

\section{Weak Properness and the Universal Curve Covering Property}
In this section we recall the definition of weak properness given in \cite{Beispielsammlung}, and 
we give another interpretation of this notion  in terms of a certain curve covering property.
We start with the curve covering property, and we first  consider algebraic prevarieties 
defined over an arbitrary algebraically closed field $\KK$. Following the 
terminology used e.g.~in \cite{Bo}, we do not require a prevariety to be irreducible. 
When we speak of a curve in a prevariety
 we mean a closed algebraic subset of pure dimension $1$. So a curve
in this sense is also not necessarily irreducible.


\begin{definition}
Let $p \colon X \to Y$ be a morphism of prevarieties.
We say that $p$ satisfies the {\it curve covering property\/} (CCP) 
if for every irreducible curve $Y'\subset Y$ and every $y\in Y'$ there
is an irreducible curve $X'\subset X$ such that $y\in p(X')\subset Y'$
and $p(X')$ is dense in $Y'$.
If the curve covering property remains true even after any
base change then we say that $p$ satisfies the {\it universal curve
covering property\/}.
\end{definition}

\begin{example}
Every surjective proper morphism of prevarieties satisfies the universal curve
covering property.
\end{example} 

\proof Let $p\colon X\to Y$ be a proper surjective morphism. 
Consider an irreducible curve $Y'$ in $Y$ and a point $y\in Y'$. Choose an
irreducible curve
$X'$ in $X$ such that $p(X')$ is dense in $Y'$. Since $p$ is proper,
the curve $X'$ must intersect the fiber of $y$. \endproof

%
%
In fact, the  curve covering property is nothing but a geometric
characterization of submersiveness:

\begin{lemma}\label{CCP = submersiv} A surjective morphism 
$p\colon X\to Y$ of  algebraic prevarieties satisfies the curve covering property  
if and only if it is submersive, i.e. if $Y$ carries
the quotient topology with respect to $p$.
\end{lemma}

\proof First suppose  that the (CCP) holds,
and consider a subset $U\subset Y$ whose preimage $p^{-1}(U)$ is open
in $X$. Let us assume that the complement $A:=Y\setminus U$ is not closed
and choose a point $y\in \overline{A}\cap U$.
Since $A=p(X\setminus p^{-1}(U))$ is constructible, there is an irreducible
curve $C_Y\subset Y$ through $y$ such that  $C_Y\cap A$ is open
and dense in $C_Y$.
Using the (CCP) we can  find an irreducible curve $C_X\subset X$ meeting
the fiber $p^{-1}(y)$ in a point $x$ such that $p(C_X)$ is dense in $C_Y$. 
That implies $x \in \overline{p^{-1}(A)}\setminus p^{-1}(A)$ contradicting
 the assumption that $p^{-1}(A)$ is closed.

Conversely suppose that $p$ is submersive and consider
an irreducible curve $C_Y$ through a point $y$ in $Y$. Since $p$ is
 surjective, the fiber of $y$ is not empty. Moreover, since
$C_Y\setminus\{y\}$ is not closed, the assumption implies that its preimage 
$p^{-1}(C_Y)\setminus p^{-1}(y)$ is also not closed. So there is a point
$x$ in the fiber of $y$ which is contained in the closure of
$p^{-1}(C_Y\setminus\{y\})$. 
This implies that there is a curve $C_X$ through $x$ in $p^{-1}(C_Y)$
intersecting the fiber of $y$ only in a finite number of
points. So the (CCP) is fulfilled.
\endproof

An example of a surjective  morphism that does not satisfy
the curve covering property is
the following:

\begin{example} Let $X$ denote the blow--up of $\KK^2$ in the origin,
and let $x=\infty$ be the point corresponding to the $e_2$--axis in the
exceptional line. Then the morphism $p\colon X\setminus\{x\}\to \KK^2$
defined by contracting the exceptional line is surjective. But there is
no curve in $X\setminus\{x\}$ covering the $e_2$--axis near the origin.
\end{example}

Here is an example of a surjective morphism with curve covering property
but such that the (CCP) does not hold universally.

\begin{example} Let $X$ denote the simple nodal curve in $\KK^2$ defined
by the equation $y^2=x^2(x+1)$. Its normalization is given by
$\nu\colon \KK^1\to X$, $t\mapsto (t^2-1,t(t^2-1))$. The  map 
$p\colon \KK^1\setminus\{-1\} \to X$ defined by $\nu$ is 
surjective and the (CCP) holds. But base change of $p$ via $\nu$
leads to a map that does not have the (CCP):

The fiber product of $\KK^1$ and $\KK^1\setminus\{-1\}$ over $X$ 
is the reducible subvariety $Y:=\{(t,t) ; t\in\KK, t\ne -1\}\cup (-1,1)$ of
$\KK^2$, and for the projection $p_1\colon Y\to \KK^1$ onto the first factor
the (CCP) does not hold.
\end{example}

From now on let us assume that all prevarieties are defined over $\CC$. 
In this case the universal curve covering property 
has a local interpretation  in terms of holomorphic germs of curves, and
for this purpose we recall some definitions from
\cite{Beispielsammlung}.  A {\it local curve} in $x \in X$ is defined
to be a holomorphic mapping germ $\gamma \colon \CC_{0} \to X'_{x}$ 
where $X'$ is an algebraic curve in $X$ through $x$. Let
$p \colon X \to Y$ be a regular map of prevarieties. We say that a
local curve $\t{\gamma} \colon \CC_{0} \to X'_{x}$ in $x \in X$ is a 
{\it weak $p$--lifting} of
a local curve $\gamma \colon \CC_{0} \to Y'_{y}$ in $y \in Y$ (where
$Y'\subset Y$ is a curve through $y$) if there is a non--constant
holomorphic mapping germ $\alpha \colon \CC_{0} \to \CC_{0}$ and a
commutative diagram
$$
\xymatrix{ \CC_{0} \ar[r]^{\t{\gamma}} \ar[d]_{\alpha} & 
X'_{x} \ar[d]^{p} \cr
 \CC_{0} \ar[r]^{\gamma} & Y'_{y} }$$
The map $p$ is called {\it
  weakly proper}, if any local curve in $Y$ admits a weak
$p$--lifting. 
A similar  notion in the context of algebraic spaces was
introduced by Koll\`ar, the so-called weak lifting property
for discrete valuation rings (see \cite{Ko},
Section~3).


%


\begin{proposition} 
For a surjective morphism $p\colon X\to Y$
of 
 complex algebraic prevarieties the following
conditions are equivalent:
\begin{enumerate}\setlength{\itemsep}{0pt}
\item $p$ is weakly proper.
\item $p$ has the universal curve covering property.
\item $p$ is submersive, and this property is preserved by every base change.
\end{enumerate} 
\end{proposition}

\proof First note that surjectivity is preserved under base change.
So the equivalence of the last two statements follows from
the Lemma~\ref{CCP = submersiv}.
Moreover, note that weak properness is preserved under base change
and clearly implies the (CCP). That shows that (i) implies (ii).

Now let us assume (ii) and conclude (i). Without loss of generality we can also assume
that $Y$ is separated. Let $\gamma$ be a local curve through a point
$y\in Y$, let $C_Y$ denote the Zariski closure of the image of $\gamma$
in $Y$. The local curve $\gamma$ factors through
the normalization $\t{C_Y}$, and after a base change we can 
achieve that $Y$ is normal and $1$-dimensional. 

By assumption there is a point $x$ in the fiber of $y$ and
an irreducible curve $C_X$ through  $x$ such that
$p(C_X)$ is dense in $Y$. Consider the normalization
$\t{C}$ of $C_X$, and choose a point $\t{x}\in\t{C}$ above $x$. 
We have an induced dominant morphism from $\t{C}$ to $Y$ mapping
$\t{x}$ to $y$.

Now we can argue locally in the analytic category. 
The germs $\t{C}_{\t{x}}$ and $Y_y$ are smooth and hence isomorphic 
to $\CC_0$. The holomorphic germ of the morphism $p$ looks like the
germ $\CC_0\to\CC_0$ defined by $z\to z^n$ for some $n\in\NN$, and
similarly the germ of $\gamma$ is of the form $\CC_0\to \CC_0$,
$z\to z^m$ for some $m$. Since both germs commute, one can choose
$\alpha=p$ and $\t{\gamma}=\gamma$ to obtain
the desired commutative diagram. So in fact (i) holds.
 \endproof

As mentioned in the introduction we want to further investigate
weakly proper toric quotients.
The main task will be to check the defining property of a categorical
quotient. In this context the following factorization result is
particularly useful
(see Proposition 1.1 in \cite{Beispielsammlung}):

\begin{proposition}\label{factmor}
  Let $p \colon X \to Y$ be a weakly proper  morphism of
  prevarieties, 
assume that $Y$ is normal and let $f \colon X \to Z$
  be a morphism into a variety $Z$. If $f$ is constant on the fibres of
  $p$, then there is a unique morphism $\t{f} \colon Y \to Z$ such
  that $f = \t{f} \circ p$. \quad$\kasten$
\end{proposition}

Note that for open surjections the above result is well--known (see
e.g. \cite{Bo}, II.6.2). 
However, 
the statement is not true in general for arbitrary morphisms~$p$.


\section{Toric Morphisms and Weak Properness}
We now come to the toric setting. Since for the proof of our
main result we need the non-separated analogues of toric varieties, 
here we first briefly summarize the basic facts about
toric prevarieties. 
Then we state the characterization
of weak properness in terms of fans.

By definition a {\it toric
prevariety\/} is a normal prevariety together with an effective action
of an algebraic torus having a dense orbit. We also fix an embedding
of the torus $T$ in the toric prevariety $X$ and denote the point in $X$
corresponding to the identity element by $x_0$.
A morphism $f\colon X\to X'$ of toric prevarieties with tori $T$ and
$T'$ respectively is called a toric morphism if $f$ maps $T$ into $T'$
and is equivariant with respect to the actions of $T$ and $T'$ respectively.
In particular, the restriction map $f|_T \colon T\to T'$ is a group
homomorphism, and we will refer to its kernel as  the {\it kernel of $f$}
and denote it by $\ker(f)$.

As in the separated case, one can associate to each toric
prevariety a convex--geo\-me\-tri\-cal object. More precisely, the category
of toric prevarieties is equivalent to the category of affine systems
of fans (see \cite{Prevarieties}). Let us recall the basic definitions.
A {\it system of fans} in a lattice $N$ is a finite collection 
$\S=(\Delta_{ij})_{i,j\in I}$ of fans $\Delta_{ij}$ in the lattice $N$
with
$$\Delta_{ij}=\Delta_{ji} \quad\text{and}\quad 
\Delta_{ij}\cap\Delta_{jk} \subset \Delta_{ik}$$
for all $i,j,k$.
Such a system of fans is called {\it affine} if for every $i\in I$
the fan $\Delta_{ii}$ consists of the faces of a single cone $\sigma_{ii}$
in $N$. 

Given an affine system of fans $\S$ in a lattice $N$, one can construct a toric
prevariety $X_{\S}$ with torus $T=N\otimes_{\ZZ} \CC^*$ by taking 
the affine toric varieties $X_i$ associated to the cones $\sigma_{ii}$
in the lattice $N$,
and glueing $X_i$ and $X_j$ 
along the open toric subvariety corresponding to the common subfan
 $\Delta_{ij}$ of $\Delta_{ii}$ and $\Delta_{jj}$ for every $i,j\in I$.
The subfans
$\Delta_{ij}$ induce a glueing relation on the set
$$\F(\S):=\{(\tau,i) ; i\in I,\tau\prec\sigma_{ii}\}$$
of labelled faces of the maximal cones occuring in $\S$, namely
$(\sigma,i)\sim(\tau,j)$ if and only if $\sigma=\tau\in\Delta_{ij}$.
There is a $1$--$1$--correspondence between the set of equivalence classes 
$\Omega(\S):=\F(\S)/\sim$ and the set of $T$--orbits in $X_{\S}$.
More precisely, every labelled face $(\tau,i)\in\F(\S)$ defines a 
distinguished point $x_{\tau}$ in the toric variety $X_i$, and
this point is identified with $x_{\tau}$ in $X_j$ if and only if $\tau\in\Delta_{ij}$.
Therefore every equivalence class $[\tau,i]\in\Omega(\S)$ defines a
distinguished point $x_{[\tau,i]}$ in $X_{\S}$, and the orbits
$T\mal x_{[\tau,i]}$, $[\tau,i]\in\Omega(\S)$, form a partition of $X_{\S}$.
The orbit structure is reflected by the partial ordering on $\Omega(\S)$
given by $[\tau,j]\prec[\sigma,i]$ if and only if $\tau\prec\sigma$ and 
$[\tau,j]=[\tau,i]$. We have
$$ T \mal x_{[\sigma,i]} \subset \b{T \mal x_{[\tau,j]}}
\Leftrightarrow [\tau,j] \prec [\sigma,i].$$

We will also need the description of toric morphisms in terms of
systems of fans. For our purposes however, it will suffice to consider
toric  morphisms from  prevarieties to  varieties. So let $X_{\S}$
be the toric prevariety arising from an affine system of fans $\S$ in a lattice
$N$ and let $X_{\Delta}$ denote the toric variety associated to
a fan $\Delta$ in a lattice $N'$. Set
$$\C(\S):=\bigcup_{i\in I} \Delta_{ii}\,.$$
Then any toric morphism
$f\colon X_{\S}\to X_{\Delta}$ corresponds to a lattice homomorphism
$F\colon N\to N'$ with the property that for every $\tau\in\C(\S)$ there is
a cone $\sigma\in\Delta$ with $F_{\RR}(\tau)\subset \sigma$.
(Here $F_{\RR}$ denotes the scalar extension of $F$ to the real
vectorspaces generated by $N$ and $N'$.)
For later use we also introduce the following notation. For
any given natural number $k$,
we denote the subset of $k$--dimensional cones in $\C(\S)$ by
$\C(\S)^k$.

The {\it support\/} of the system of fans $\S$ is defined to be
$|\S|=\bigcup_{i\in I} |\Delta_{ii}|$. For toric morphisms 
weak properness can be characterized as follows
(see Proposition 1.2 in \cite{Beispielsammlung}):

\begin{proposition}\label{weakpropchar}
  A toric morphism $f\colon X_{\S}\to X_{\Delta}$ from a toric prevariety
to a toric variety is weakly proper if and only if the 
associated lattice homomorphism $F$ induces a surjection on the
supports of the corresponding systems of fans
$$F_{\RR}(\vert \S \vert) = \vert \Delta \vert\,. \quad \kasten$$
\end{proposition}

\section{Partition of the Support defined by a Morphism}

\noindent
Let us consider a morphism $f\colon X\to Z$ from a
toric prevariety $X=X_{\S}$ arising from an affine system of fans
$\S$ in a lattice $N$ to an algebraic variety $Z$
that is not necessarily toric. As we will see, such a morphism
defines an equivalence relation on the set $\C(\S)$ of cones
occuring in $\S$, and the supports of the equivalence classes
form a partition of the support of $\S$ in finitely many
subsets.
Whether or not $f$ factors through a given toric
morphism can be expressed in terms of this partition.

We will consider two elements of $\Omega(\S)$ as equivalent with
respect to $f$ if the corresponding parametrized orbits are mapped
by $f$ to the same parametrized set in $Z$, or more precisely if
the composition of $f$ with the orbit maps of the
corresponding distinguished points yield the same map on $T$.  

\begin{definition} The morphism $f$ induces an equivalence relation
on the set $\Omega(\S)$, \break
namely:
$$ [\sigma,i] \sim_f [\sigma',j] \quad \Leftrightarrow \quad
 f(t\mal x_{[\sigma,i]})=f(t\mal x_{[\sigma',j]}) 
\quad\text{for all $t\in T$.}$$
\end{definition}

Note that if the same cone appears with two different labels in $\Omega(\S)$,
$[\sigma,i]$ and $[\sigma,j]$ say, then $[\sigma,i]\sim_f [\sigma,j]$. That is
an immediate consequence of the following remark. For a cone $\sigma$,
let $\sigma^{\circ}$ denote its relative interior. We observe:

\begin{remark}\label{ueberlappende Kegel} If 
for $[\sigma,i]$, 
$[\sigma',j]\in \Omega(\S)$
we have $\sigma^{\circ}\cap (\sigma')^{\circ}\ne \emptyset$, 
then $[\sigma,i]\sim_f [\sigma',j]$.
\end{remark}

\proof Choose $v\in N\cap \sigma^{\circ}\cap (\sigma')^{\circ}$. Let
 $\lambda_{v}$ denote the corresponding
one--parameter--subgroup of $T$ and fix $t\in T$. In the affine
chart $X_i$, $\lim_{s\to 0} t\mal\lambda_v(s)\mal x_0=t\mal x_{[\sigma,i]}$
holds, whereas in $X_j$ we have
$\lim_{s\to 0} t\mal\lambda_v(s)\mal x_0=t\mal x_{[\sigma',j]}$.
Since $f$ is continuous and we assumed
$Z$ to be separated, this implies the claim.
\endproof

So in fact, $f$ induces an equivalence relation on the set $\C(\S)$,
and we will also denote this relation by $\sim_f$.
We define
the {\it support of the equivalence class} of $\sigma\in\C(\S)$ by
$$|\sigma|_f:=\bigcup_{\sigma'\sim_f \sigma} (\sigma')^{\circ} .$$
As an immediate consequence of Remark \ref{ueberlappende Kegel}
we obtain the following

\begin{remark} The subsets $|\sigma|_f$, $\sigma\in\C(\S)$, 
form a partition of the support of $\S$.\endproof
\end{remark}



\begin{example}\label{regions of tormor}
Consider a fan $\Delta$ in a lattice $N'$, and let $X_{\Delta}$ denote
the corresponding toric variety. 
Let $f\colon X_{\S} \to X_{\Delta}$ be a toric morphism,
and let $F\colon N\to N'$ denote the associated lattice homomorphism.
Then the equivalence classes of $\C(\S)$ correspond to the elements of
$\Delta$ that meet $F(|\S|)$, more precisely
$\sigma\sim_f\sigma'$ if and only if there is a cone $\tau\in\Delta$ with
$F(\sigma^{\circ})\subset \tau^{\circ} \supset F((\sigma')^{\circ})$.
The supports of the equivalence classes are the sets 
$F^{-1}(\tau^{\circ}) \cap |\S|$,
$\tau\in\Delta$. 

For example, let us look at the toric morphism $f\colon \CC^2\to \CC$
given by the projection on the first factor. The system of fans corresponding
to $\CC^2$ is the fan $\Delta:=\{\sigma,\tau_1,\tau_2,0\}$ in $\ZZ^2$,
where $\sigma:=\cone(e_1,e_2)$, $\tau_i:=\cone(e_i)$ for $i=1,2$, 
 and $\CC$ arises from the fan $\Delta':=\{\cone(e_1),0\}$ in $\ZZ$.
The lattice homomorphism corresponding to $f$
is the projection $F\colon \ZZ^2\to \ZZ$.
So in this case we have $0\sim_f \tau_2$ and $\tau_1\sim_f\sigma$.
The supports of the equivalence classes in $|\S|=\sigma$ are
$|\tau_2|_f=\tau_2$ and
$|\sigma|_f=\sigma^{\circ}\cup \tau_1^{\circ}$.
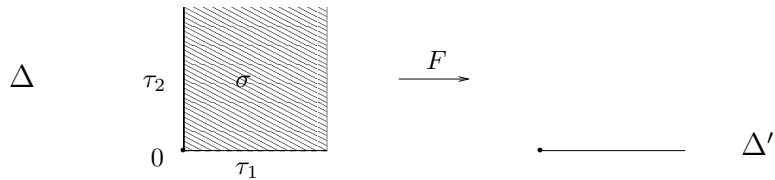
\begin{center}
\input{regions.pstex_t}  
\end{center}
\end{example}

\bigskip

For a given cone $\sigma\in \C(\S)$, let $T(|\sigma|_f)$ denote the 
subtorus of $T$ corresponding to the sublattice obtained by intersecting
$N$ with the linear subspace $\lin |\sigma|_f$ generated in $N_{\RR}$
by the set $|\sigma|_f$.
Then $T(|\sigma|_f)$ is generated by all isotropy subgroups
$T_{x_{[\sigma',i]}}$, where $[\sigma',i]\in\Omega(\S)$ and
$\sigma'\sim_f\sigma$. 

\begin{remark}\label{relstab} 
We have $f(t\mal t'\mal x_{[\sigma,i]})= f(t\mal x_{[\sigma,i]})$
for every $t'\in T(|\sigma|_f)$. In particular, $f$ is invariant
with respect to the action of $T(|0|_f)$.
\end{remark}

\proof To see this, choose cones $\sigma_1,\dots,\sigma_s$
in the $f$--equivalence class of $\sigma$ 
that are generating $\lin |\sigma|_f$
as a vector space. If $[\sigma_j,i_j]\in\Omega(\S)$,
then $f(t\mal x_{[\sigma_j,i_j]})=f(t\mal x_{[\sigma,i]})$ for every $t\in T$,
since $\sigma_j\sim_f \sigma$. So we can conclude
that $f(t'\mal t\mal x_{[\sigma,i]})=f(t\mal x_{[\sigma,i]})$ for 
every $t'$ in the stabilizer $T_j$ of the point $x_{[\sigma_j,i_j]}$.
Since the subtori $T_j$ generate $T(|\sigma|_f)$ that implies the claim.
\endproof

Let us now consider a dominating toric morphism 
$p\colon X_{\S}\to X_{\Delta}$
to some toric variety $X_{\Delta}$, and assume that the associated
lattice homomorphism $P\colon N\to N'$
is surjective. That means that the kernel  of the 
homomorphism of tori $T\to T'$ induced by $P$ is connected. We will denote
this kernel by $\ker(p)\subset T$.
Assume that $p$ is weakly proper.
 With the above notations we can describe the 
fibers of $p$  as follows:

\begin{lemma}\label{Fasern} 
Two points $x,y$ lie in the same fiber of p if and only if
there are elements $[\sigma_1,i],[\sigma_2,j]\in\Omega(\S)$ 
with $\sigma_1\sim_p\sigma_2$, $t\in T$ and $t_1\in T(|\sigma_1|_p)\mal \ker(p)$
with $x=t\mal x_{[\sigma_1,i]}$ and $y=t\mal t_1\mal x_{[\sigma_2,j]}$.
\end{lemma} 

\proof Let $p|_T \colon T \to T'$ denote the restriction homomorphism of  $p$
to the big tori of $X_{\S}$ and $X_{\Delta}$ respectively.
Any point $z\in X_{\Delta}$ is of the form $z=t'\mal x_{\sigma'}$
for some $t'\in T'$ and $\sigma'\in\Delta$.
Then the $p$-fibre of the point $z$ is
$$
p^{-1}(z)=p^{-1}(t'\cdot x_{\sigma'}) = \bigcup_{P_{\RR}(\sigma)^{\circ}
  \subset (\sigma')^{\circ}} (p|_T)^{-1}(t' \mal T'_{x_{\sigma'}})
\mal x_{[\sigma,i]}$$
(see \cite{Beispielsammlung}, Proposition 3.5). As described in Example
\ref{regions of tormor}, since $p$ is surjective, the given 
cone $\sigma'\in\Delta$ defines a $p$--equivalence class, represented
by $\sigma_1\in\C(\S)$ say, and $|\sigma_1|_p=
|\S|\cap P_{\RR}^{-1}((\sigma')^\circ)$.
We obtain:
$$
p^{-1}(t'\cdot x_{\sigma'}) = \bigcup_{\sigma\sim_p
  \sigma_1} (p|_T)^{-1}(t' \mal {T'}_{x_{\sigma'}})
\mal x_{[\sigma,i]}.$$
To prove the lemma it suffices to show that
 $$(p|_T)^{-1}(T'_{x_{\sigma'}})=T(|\sigma_1|_p)\cdot
\ker(p)\,.$$
First note that the subtorus $T'_{x_{\sigma'}}$ corresponds
to the sublattice of $N'$ defined by $\lin\sigma'$. 
By Proposition \ref{weakpropchar}, since $p$ is weakly proper
 we have $P_{\RR}(|\S|)=|\Delta|$.
That implies  $P_{\RR}(|\sigma_1|_p)=(\sigma')^{\circ}$ and hence
$P_{\RR}(\lin {|\sigma_1|_p})=\lin \sigma'$. Therefore
$P^{-1}_{\RR}(\lin\sigma')=\lin |\sigma_1|_p + \ker(P_{\RR})$.
Since we assumed $P$ to be surjective, that proves the claim.
\endproof

From this description of the fibers of a weakly proper toric morphism 
we obtain the following factorization criterion. 

\begin{lemma}\label{factorization criterion} Let $f\colon X_{\S}\to Z$ be
a morphism to a variety $Z$, and let
$p\colon X_{\S}\to Y$ be a  toric morphism with the universal curve
covering property to a toric
variety $Y$, such that $\ker(p)$ is connected. 
Then the following statements
are equivalent:
\begin{enumerate}
\item The morphism $f$ factors through $p$.
\item $f(t_1\cdot x_{[\sigma_1,i]})=f(t_2\cdot x_{[\sigma_2,j]})$ for all
$[\sigma_1,i],[\sigma_2,j]\in\Omega(\S)$ with $\sigma_1\sim_p \sigma_2$ and
 $t_1,t_2\in T$ with $t_2^{-1} t_1\in T(|\sigma_1|_p)\mal\ker(p)$.
\item $f$ is  $\ker(p)$--invariant, and for 
$\sigma,\sigma'\in\C(\S)$, whenever
$\sigma\sim_p \sigma'$ then
$\sigma\sim_f \sigma'$. 
\end{enumerate}
\end{lemma}

\proof By Proposition \ref{factmor}, $f$ factors through $p$ if and only if
$f$ is constant on the fibers of $p$. 
From the above description of the fibers of $p$ it follows immediately that
 (i) and (ii) are equivalent.
Now let us  assume that (i) holds. Then in
particular $f$ is $\ker(p)$--invariant. Moreover, it can be read off directly
from (ii) that
$\sigma\sim_f\sigma'$ whenever $\sigma\sim_p\sigma'$.  

Conversely, assume that (iii) holds, and let $\sigma\in\C(\S)$.
Then the $f$--equivalence class of $\sigma$ contains the $p$--equivalence
class of $\sigma$, and hence $T(|\sigma|_p)\subset T(|\sigma|_f)$. 
Now (ii) follows from Remark \ref{relstab}.
\endproof

We want to apply this lemma to the following situation. 
Let $\S$ be an affine system of fans in a lattice
$N$ and let $\sigma$
be a  convex cone in $N$, not necessarily contained
in $\C(\S)$ but with $\sigma\subset |\S|$. Let $\sigma_0$ denote
the minimal face of $\sigma$. Then $\sigma_0$ is a linear subspace
of $N_{\RR}$, and the sublattice $\sigma_0\cap N$ defines a subtorus
$H$ of the big torus of $X_{\S}$. 
 
Define a system of fans $\S\cap\sigma=(\Delta'_{ij})_{i,j\in I}$, 
in $N$ by setting $\sigma_{ii}':=\sigma_{ii}\cap \sigma$ and
$\Delta_{ij}':=\{\tau\cap\sigma ; \tau\in\Delta_{ij}\}$.
Then  the identity homomorphism $\id_N$  defines a toric morphism
$q\colon X_{\S\cap\sigma}\to X_{\S}$, and the projection 
$P\colon N\to N/(\sigma_0\cap N)$ defines a toric morphism
$p\colon X_{\S\cap\sigma} \to X_{P(\sigma)}$ to the affine 
toric variety $X_{P(\sigma)}$ associated to the cone $P(\sigma)$ in 
$N/(\sigma_0\cap N)$.

\begin{corollary}\label{indMor} 
Let $f\colon X_{\S}\to Z$ be a morphism to some variety $Z$, and assume
that for every face $\tau$ of $\sigma$, the set $\tau^{\circ}$ is contained in
the support of an $f$--equivalence class of $\C(\S)$.
Then the morphism $f$ is $H $--invariant, and
 there is a unique  morphism $f_{\sigma}\colon X_{P(\sigma)}\to Z$ 
such that the following diagram is commutative:
$$\xymatrix{
X_{\S\cap\sigma} \ar[r]^{q} \ar[dr]_{p} & X_{\S}  
\ar[r]^{f} & \ Z \cr
& X_{P(\sigma)} \ar[ur]_{f_{\sigma}} & \cr
}
\qquad.$$
\end{corollary}

\proof The minimal face $\sigma_0$ is a linear subspace and
therefore it coincides with its relative interior.
By assumption, $\sigma_0$ is contained in the support of an
$f$--equivalence class, and that means, that
$\sigma_0\subset |0|_f$. Therefore $H\subset T(|0|_f)$, 
and Remark \ref{relstab} implies that $f$ is
$H$--invariant.

By construction, $P$ is surjective and weakly proper.
So we can apply Lemma \ref{factorization criterion}, and conclude 
that the morphism $f\circ q$ factors uniquely
through $p$. That means that there is a morphism 
$f_{\sigma}\colon X_{P(\sigma)}\to Z$ with $f\circ q=f_{\sigma}\circ p$.
\endproof

\section{Closures of Cones in the Support}

\noindent 
The following two sections contain preparations
 for the proof of the main result given in Section 6.
First we  analyze morphisms from a toric variety
corresponding to a fan with  
just two maximal cones in a special position.

\begin{lemma}\label{zwei Kegel} 
  Let $\Delta$ be a fan with only two maximal cones $\sigma_1$ and $\sigma_2$.
  Assume that there are faces $\tau_{i}$ of $\sigma_i$
   and vectors $v\in(\tau_1\cap\tau_2)^{\circ}\cap N$,
$v+v'\in
  \tau_1^{\circ} \cap N$, $v-v'\in\tau_2^{\circ}\cap N$,
  $w\in(\sigma_{1}\cap\sigma_{2})^{\circ}\cap N$ with
  $v+v'+w\in\sigma_{1}^{\circ}$. 
  Let $f\colon X_{\Delta}\to Z$ be a morphism such that
$\sigma_{1}\sim_f(\sigma_1\cap\sigma_{2})$.
  Then $\tau_1\sim_f(\tau_1\cap\tau_2)$.
\end{lemma}

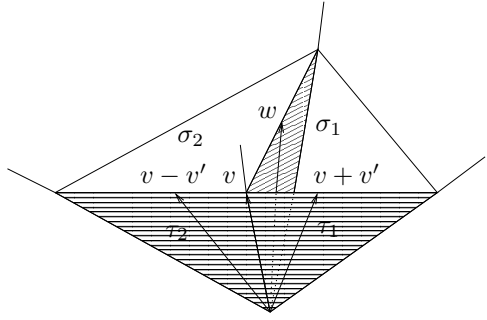
\begin{center}
\input{zweikegel.pstex_t}
\end{center}

\proof Set $\rho:=\tau_1\cap\tau_2$.
We have to show that for every $t\in T$ the following holds:
 $$f(t\mal x_{\rho})=
   f(t\mal x_{\tau_1})
\;.$$
For a given $t\in T$, consider the morphism $f_t\colon X_{\Delta}\to Z$,
defined by $f_t(x):=f(t\mal x)$. Then clearly $f_t$ satisfies the
same assumption as $f$, i.e.
$\sigma_{1}\sim_{f_t}(\sigma_1\cap\sigma_{2})$. Therefore it suffices to show
the claim for $t=1$.

Let $V$ denote the
closure of the orbit $T\mal x_{\rho}$ in $X_{\Delta}$.
Then $V$ is a toric variety with respect to the torus $T/T_{x_{\rho}}$,
and it corresponds to the fan obtained by projecting the star of
$\rho$ in $\Delta$ to $N/(\lin\rho\cap N)$.
The lattice homomorphism
$F\colon \ZZ^2\to N/(\lin(\rho)\cap N)$,
defined by $F(e_1)=\b{v'}$ and $F(e_2)=\b{w}$, defines a toric morphism 
$\varphi\colon \PP_1\times\CC\to V$. This toric morphism
has the following properties:
$$\varphi([r_0,r_1],s)= \cases{ \lambda_{v'}(r_1/r_0) \mal \lambda_{w}(s)
  \mal x_{\rho} & if $r_0,r_1,s \ne 0$, \cr 
\lambda_{w}(s) \mal x_{\tau_1} & if $r_1= 0$, $s\ne 0$, \cr
  \lambda_{v'}(r_1/r_0)\mal x_{\sigma_1\cap\sigma_2} & 
   if $r_0,r_1\ne0$, $s=0$ \cr }\,.$$

Now consider the composition $\psi=f\circ\varphi\colon \PP_1\times\CC\to Z$.
By definition, $v'\in\lin \sigma_1$ and hence 
$\lambda_{v'}(\CC^*)\subset T(|\sigma_1|)$.
Since we assumed  $\sigma_1\sim_f(\sigma_1\cap\sigma_2)$, this implies 
$f(\lambda_{v'}(r)\mal x_{\sigma_1\cap\sigma_2})=f(x_{\sigma_1})$ for all $r\in\CC^*$.
So $\psi(\PP_1\times\{0\})=f(x_{\sigma_1})$, i.e.~$\psi$ contracts
the curve $\PP_1\times\{0\}$ to a point.
With the following Lemma \ref{Kurvenschar} we can conclude
that $\psi$ in fact does not depend on the first coordinate,
and for $s=1$ we obtain that $f(x_{\rho})=f(x_{\tau_1})$.
\endproof

In the above proof we used a general fact about 
morphisms from $\PP_1\times\CC$.

\begin{lemma}\label{Kurvenschar} Let 
$\psi\colon \PP_1 \times \CC \to Z$ be a morphism to a variety $Z$ with
$\psi(\PP_1\times\{0\})=z$ for some $z\in Z$. Then $\psi$ is  constant
on $\PP_1\times \{s\}$ for every $s\in\CC$, i.e.
$\psi$ does not depend on the first coordinate. 
\end{lemma}

\proof Choose an open affine neighbourhood $W$ of $z$ in $Z$ and set $Y
:= (\PP_1 \times \CC) \setminus \psi^{-1}(W)$. Consider the
projection $\pr \colon \PP_1 \times \CC  \to \CC$.
Since $\PP_1$ is complete, $\pr(Y)$ is closed in $\CC$.
Moreover, we have $0 \not\in \pr(Y)$. So
$W_0 := \CC \setminus \pr(Y)$ is an open neighbourhood of $0$ in $\CC$,
and by definition
$$
 \PP_1 \times W_0 \;=\; \pr^{-1}(W_0) \;
  \subset \; \psi^{-1}(W). $$
So  by restriction we obtain a morphism 
$\psi\colon \PP_1\times W_0 \to W$.
Since we chose $W$ to be affine, $\psi$ maps $\PP_1 \times \{s\}$ to a
point for every $s \in W_{0}$. So for continuity reasons, $\psi$
does not depend on the first coordinate.
  \endproof

\goodbreak

Now we consider a system of fans $\S$ with convex support,
and a morphism $f\colon X_{\S}\to Z$
to some variety $Z$. We apply Lemma \ref{zwei Kegel} to prove 
the following:

\begin{proposition}\label{Kegelabschluss}
Let $\sigma\subset |\S|$ be a rational (not necessarily strictly) convex cone. 
Suppose that $\sigma^{\circ}$ is contained in the support of an 
$f$--equivalence class of $\C(\S)$.  Then  for every face $\tau$ of $\sigma$,
the relative interior $\tau^{\circ}$ is also contained in the support
of an $f$--equivalence class.
\end{proposition}

\proof 
By induction on $n:=\dim \sigma$ we will show that the assertion is
true for all one--codimensional faces of $\sigma$. If $\sigma$ 
is one--dimensional
there is nothing to show. So assume  that $n\geq 2$,
and let $\tau$ be a face of $\sigma$ of dimension $n-1$.
W.l.o.g. we can assume that $\dim\sigma_{ii}=\dim\sigma$ for all $i$
since the cones of maximal dimension cover $\sigma$. We reduce
the induction step to proving the following 
\begin{description}
\item{{\bf Claim:}\enspace} For every  cone $\tau_1\in\C(\S\cap\tau)^{n-1}$
we have $\tau_1\cap\tau^{\circ}\subset |\tau_1|_f \,.$
\end{description}

From this claim it follows that
$|\tau_1|_f\cap \tau^{\circ}$ is relatively closed in $\tau^{\circ}$.
The $(n-1)$--dimensional cones  in $\C(\S\cap\tau)$
cover $\tau$, and we obtain a partition of $\tau^{\circ}$ in relatively
closed subsets of the form  $|\tau_1|_f\cap \tau^{\circ}$, where
$\tau_1\in \C(\S\cap\tau)^{n-1}$. Since $\tau^{\circ}$ is connected
this implies that one of the subsets  actually equals  $\tau^{\circ}$,
and that means that $\tau^{\circ}$ is contained in the support of
an $f$--equivalence class.

So to prove the proposition it suffices to show the claim.
Let us assume the claim were not true. Then since all the cones in question
are rational, we can find a rational vector
$v$ in the boundary of $\tau_1$ with 
$v\in\tau^{\circ}\setminus |\tau_1|_f$.
Note that $v$ may be
zero, if the cone $\tau$ is not strictly convex. 

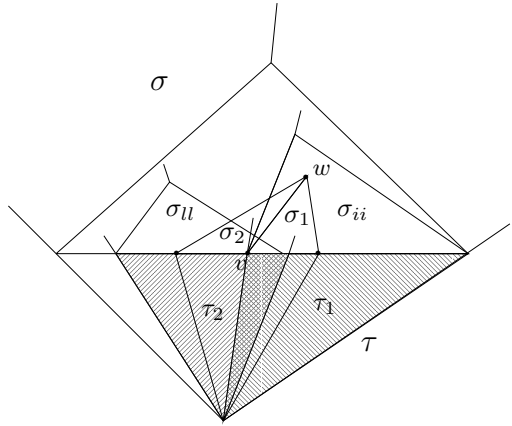
\begin{center}
\input{convex.pstex_t}
\end{center}
\goodbreak

Now choose a vector $v'\in \lin\tau$ such that
 $v+v'\in\tau_1^{\circ}$.
For sufficiently small $\epsilon$ the ball $B_{\epsilon}(v)$ 
of radius $\epsilon$ in $\lin\tau$ around $v$ is already covered by
those  cones in $\C(\S\cap\tau)^{n-1}$  that contain $v$.
Therefore if we choose $v'$ of length $\leq \epsilon$, then
$v-v'\in \tau_2$ for some $\tau_2\in \C(\S\cap\tau)^{n-1}$ containing $v$.

Moreover, for $i=1,2$ we can construct  $(n-1)$--dimensional
rational cones $\tau_i'\subset\tau_i$ such that $v+v'\in(\tau_1')^{\circ}$,
$v-v'\in(\tau_2')^{\circ}$ and
$\tau_1'\cap\tau_2'=\cone(v)$ is a common face of $\tau_1'$ and
$\tau_2'$. To do so, we can first choose a hyperplane through $v$
separating $v+v'$ and $v-v'$, and then choose appropriate simplicial
rational cones around $v+v'$ and $v-v'$ respectively that lie entirely on one
side of the hyperplane, and then form the convex hull with $v$.

By assumption there are $n$--dimensional cones $\sigma_{ii}$ and 
$\sigma_{ll}$ having $\tau_1$
and $\tau_2$ respectively as a face. 
Now choose $w\in \sigma_{ii}^{\circ}\cap N$,
and set $\sigma_1:=\cone(\tau_1', w)$ and $\sigma_2:=\cone(\tau_2',w)$. 
Let $X_{\sigma_i}$ denote the toric variety associated to $\sigma_i$
in $N$ ($i=1,2$). Since the cone $\sigma_1$ is 
contained in the cone $\sigma_{ii}$, the identity on $N$ defines a toric
morphism from  $X_{\sigma_1}$ to the affine chart $X_i:=X_{\sigma_{ii}}$
of $X_{\S}$.
The composition with $f$ yields a morphism
$f_1$ from $X_{\sigma_1}$ to $Z$.

The cone  $\sigma_2$
has  the following properties: 
$\sigma_2\cap\tau=\tau_2'\subset \tau_j$ and 
$\sigma_2\backslash \tau\subset\sigma^{\circ}$. 
Therefore the relative interior of every face of $\sigma_2$ 
is contained in an $f$--equivalence class. Note that $\sigma_2$ is strictly
convex, since by construction $v'\notin \lin(v)$.  Hence by 
Corollary \ref{indMor},
 $f$ defines a morphism $f_2$ from $X_{\sigma_2}$ to $Z$ such that
the following diagram is commutative:
$$\xymatrix{X_{\S\cap\sigma_2} \ar[r]  \ar[dr] & X_{\S}  
\ar[r]^{f} & \ Z \cr
& X_{\sigma_2} \ar[ur]_{f_2}& \cr
}\qquad.$$

Now consider the toric variety $X_{\Delta}$ corresponding to the fan with
the two maximal cones $\sigma_1$ and $\sigma_2$ in $N$. It follows from 
the above commutative diagram that the two morphisms 
$f_1$, $f_2$ coincide on the intersection  of the two maximal
affine charts, and so they glue together to a 
morphism $f'\colon X_{\Delta}\to Z$.
Since $\sigma_1^{\circ}, \sigma_2^{\circ}\subset \sigma^{\circ}$,
by assumption we have  $\sigma_1\sim_{f'}(\sigma_1\cap\sigma_2)$.  

The cones $\sigma_i$ with the faces $\tau_i'$ satisfy
all the conditions of Lemma \ref{zwei Kegel}, and we obtain
$\tau_1'\sim_{f'}\rho$, where $\rho:=\tau_1'\cap\tau_2'=\cone(v)$
and hence $v\in |\tau_1'|_{f'}$.
Since $(\tau_1')^{\circ}\cap \tau_1^{\circ}\ne \emptyset$,
it follows that  $v\in |\tau_1|_f$,
which is a contradiction.
\endproof

\section{The Convex Hull of Two Cones}

\noindent
As in the previous section we consider a system of fans
with convex support $|\S|$ in a lattice $N$. 
Let $\sigma,\tau\in\C(\S)$ be two cones with
$\sigma\cap \tau^{\circ}\ne\emptyset$.
Then since we assumed $|\S|$ to be convex, we have
$\sigma+\tau\subset |\S|$. Moreover, it follows
that $\sigma^{\circ}\subset (\sigma+\tau)^{\circ}=
\sigma^{\circ}+\tau^{\circ}$.

A first observation is the following:

\begin{lemma}\label{Nachbarkegel}  Let
$f\colon X_{\S}\to Z$ be a morphism from the  toric prevariety $X_{\S}$
associated to $\S$ to some variety $Z$. Suppose that
 $\sigma_1$ is a cone in $\C(\S)$ with 
$\sigma_1\cap\tau^{\circ}\ne\emptyset$ and
$(\sigma_1)^{\circ}\subset \sigma^{\circ}+\tau^{\circ}$. Then 
$\sigma_1\sim_f \sigma$.
\end{lemma}

\proof The assumptions on $\sigma_1$ imply that there is a face $\tau_1$ of
$\sigma_1$ with $\tau_1^{\circ}\cap \tau^{\circ}\ne\emptyset$.
Similarly, there is a face $\tau'$ of $\sigma$ with
$(\tau')^{\circ}\cap\tau^{\circ}\ne\emptyset$. Therefore
$\tau_1\sim_f\tau\sim_f\tau'$. 
Now choose  $v_1\in \sigma_1^{\circ}\cap (\sigma^{\circ}+\tau^{\circ})$,
and write $v_1$ in the form $v_1=w+v'$, where $w\in\sigma^{\circ}$, $v'\in\tau^{\circ}$.

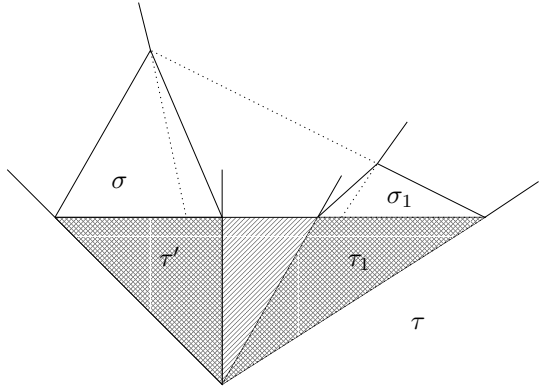
\begin{center}
\input{convhull.pstex_t}
\end{center}

Suppose that
$[\sigma,i],[\tau,j],[\sigma_1,k]\in\Omega(\S)$ and
fix $t\in T$. Since $\tau_1\sim_f\tau\sim_f\tau'$,  for all $s\in\CC^*$
the following holds:
$$
f(t\mal\lambda_{v_1}(s)\mal  x_{[\tau_1,k]})=
f(t\mal\lambda_{v_1}(s)\mal  x_{[\tau,j]})=
f(t\mal\lambda_{w}(s)\mal  x_{[\tau,j]})=
f(t\mal\lambda_{w}(s)\mal  x_{[\tau',i]})\,.$$

In the affine chart
$X_{[\sigma_1,k]}$ we have 
$\lim_{s\to\infty} t\mal\lambda_{v_1}(s)\mal x_{[\tau_1,k]}= 
t\mal x_{[\sigma_1,k]}$.
On the other hand, in the affine chart $X_{[\sigma,i]}$ we have 
$\lim_{s\to\infty} t\mal\lambda_{w}(s)\mal x_{[\tau',i]}= 
t\mal x_{[\sigma,i]}$.
That implies $f(t\mal  x_{[\sigma_1,k]})=
f(t\mal x_{[\sigma,i]})$ for all $t$, and hence $\sigma\sim_f\sigma_1$.
\endproof

\begin{proposition}\label{Kegelvergroesserung}
Let $\S=(\Delta_{ij})_{i,j\in I}$ be an affine system of 
fans with convex support, and
let $\sigma,\tau$ be cones in $\C(\S)$ with
$\sigma\cap \tau^{\circ}\ne\emptyset$. Then one of the
following two assertions holds:
\begin{enumerate}
\item
There is a nonzero linear subspace
$L\subset (\sigma+\tau)$  and a cone $\sigma_1\in \C(\S)$ such that
$L\subset |\sigma_1|_f$ for every morphism
$f\colon X_{\S}\to Z$ to a variety $Z$, or 
\item
$(\sigma+\tau)^{\circ}\subset |\sigma|_f$ for
every morphism
$f\colon X_{\S}\to Z$ to a variety $Z$.
\end{enumerate}
\end{proposition}

\proof 
We will prove the proposition by induction on the dimension 
of the support $n:=\dim\lin|\S|$ and
the number $|\C(\S)^{n}|$ of $n$--dimensional cones in $\C(\S)$.
Let us suppose that the first assertion does not hold. 
Note that this implies that the analogous assertion also does not
hold for any system of the form $\S\cap\delta$ obtained by intersecting
$\S$ with a cone $\delta\subset |\S|$.

So without loss of generality we may  assume that
$|\S|=\sigma+\tau$. Moreover, it suffices to show assertion (ii)
under the additional assumption
\begin{equation}
\sigma+\tau=\rho + \tau, \quad\text{for some ray $\rho$ of $\sigma$}.
\end{equation} 

To see this, suppose that (ii) is true whenever the extra condition (1)
holds. Choose a ray  $\rho\in \sigma\setminus\tau$ and consider the
cone $\sigma':=\rho+(\sigma\cap\tau)$. Our assumption implies
that $(\sigma'+\tau)^{\circ}\subset |\sigma|_f$ for every morphism
$f\colon X_{\S}\to Z$ to a variety $Z$. And if $v\in\sigma\cap\tau^{\circ}$,
then any point $w\in \cone(\rho,v)^{\circ}\subset \sigma'$ is contained
in the relative interior of $\tau+\rho$. Therefore we can replace
$\tau$ in $\C(S)$ by $\tau+\rho$, glued to any other cone along the
origin only. By recursion over the rays we obtain the claim.

Now we further reduce the situation to the special case that
\begin{equation}
\dim(\sigma\cap\tau)=\dim\tau.
\end{equation}

Suppose that assertion (ii) always holds if (2) is true. Consider
a  pair of cones $\sigma,\tau$ as in the proposition and satisfying 
condition (1). Then $\dim\tau=n$ or $n-1$, where
$n=\dim|\S|$. Choose a point $v_1\in\sigma^{\circ}$.
Since $|\S|=\sigma+\tau$, we can find an $n$-dimensional cone 
$\sigma_1\in\C(\S)$ containing $v_1$ such that condition (2) holds for 
$\sigma_1$ and $\sigma$. So using (ii)
we may  replace $\sigma$ in $\S$ by $\sigma_1+\sigma$. In other words,
we can assume that $\dim\sigma=n$.

If also $\tau$ is $n$-dimensional we obtain 
$\tau^{\circ}\cap\sigma^{\circ}\ne \emptyset$, and that implies
condition (2) for $\sigma$ and $\tau$.
Otherwise $\tau$ must be a facet of $|\S|$.
In that case consider a point $v\in\sigma\cap\tau^{\circ}$. As above
 we can find an $n$-dimensional cone $\tau_1\in\C(\S)$ containing $v$
such that condition (2) holds for $\tau_1$ and $\tau$. So using (ii)
we may assume that $\tau_1$ contains $\tau$ as a facet. 
From the fact that $\tau_1\subset\rho+\tau$ we can conclude that 
$\cone(v,\rho)$ meets $\tau_1$ in its relative interior. That implies
$\sigma\cap\tau_1^{\circ}\ne\emptyset$, and we can replace $\tau$
by $\tau_1$. Since $\dim\tau_1=n$, condition (2) follows as before.

From now on let us  assume that conditions (1) and (2) are satisfied.
So we are left with two possibilities: either $\dim\tau=\dim\sigma$
or $\dim\tau=\dim\sigma-1$.

Let us first consider the case $\dim\tau=\dim\sigma$.
If $\sigma\cup\tau$ is convex then there is nothing to show, since
then we have $(\sigma\cup\tau)^{\circ}=
\sigma^{\circ}\cup\tau^{\circ}\subset |\sigma|_f$.

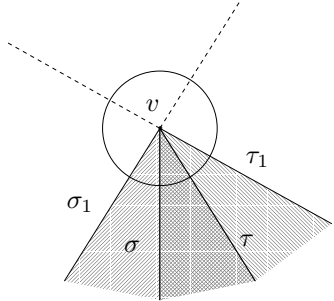
\begin{center}
\input{konvexerFall.pstex_t}
\end{center}

\noindent
(For consider $v\in\partial\sigma\cap\partial\tau$.
Then there are facets $\sigma_1\prec\sigma$ and $\tau_1\prec\tau$
with $v\in\sigma_1\cap\tau_1$.  Choose defining hyperplanes $u\in\sigma^{\vee}$
and $w\in\tau^{\vee}$ with $u^{\perp}\cap\sigma=\sigma_1$
and $w^{\perp}\cap\tau=\tau_1$. Here $\sigma^{\vee}$ and
$\tau^{\vee}$ denote the dual cones of $\sigma$ and $\tau$ respectively.
In a ball of sufficiently small radius around $v$ in $\lin\sigma$ we find
a point $v'$ with $u(v')<0$ and $w(v')<0$ and hence $v'\notin \tau\cup\sigma$.
This shows that $v$ cannot lie in the relative interior of $\sigma\cup\tau$.)

So let us  assume that $\sigma\cup\tau$ is not convex. Then there is
a facet $\sigma_1\prec\sigma$, defined by a hyperplane $u\in\sigma^{\vee}$
such that $u^{\perp}\cap \tau^{\circ}\ne\emptyset$ and 
$\sigma_1\not\subset \tau$.

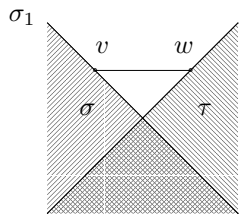
\begin{center}
\input{konkaveEcke.pstex_t}
\end{center}

\noindent
(To see this, choose a point $v\in\partial\sigma$, $w\in\partial\tau$ 
such that the segment $[v,w]$ joining $v$ and $w$ intersects $\sigma\cup\tau$
only in $\{v,w\}$. The point $v$ lies on a facet $\sigma_1\prec\sigma$,
and since $w\notin\sigma$, we can choose a defining
 hyperplane $u\in\sigma^{\vee}$ of $\sigma_1$ such that $u(w)<0$.
On the other hand, $u(\tau^{\circ}\cap\sigma^{\circ})>0$,
and so $u^{\perp}\cap\tau^{\circ}\ne\emptyset$.)

Now we decompose $\sigma+\tau$ along the defining hyperplane $u^{\perp}$
of $\sigma_1$.
Since $\sigma_1\not\subset\tau$ and $\sigma=\rho+(\sigma\cap\tau)$, 
we have $\rho\subset\sigma_1$ and therefore
$\sigma+\tau=\sigma_1+\tau$.
Set 
$$\tau_0:=\tau\cap u^{\perp}, \quad
\tau_1:=\{v\in \tau; u(v)\geq 0\}\quad \text{and}\quad
\tau_2:=\{v\in \tau; u(v)\leq 0\}\,.$$
Then $\tau_1$ and $\tau_2$ are cones that intersect in the common face
$\tau_0$. Moreover, since $\sigma+\tau=\sigma_1+\tau$,
we have $\sigma+\tau=(\sigma_1+\tau_1)\cup(\sigma_1+\tau_2)=
(\sigma+\tau_1)\cup (\sigma_1+\tau_2)$.

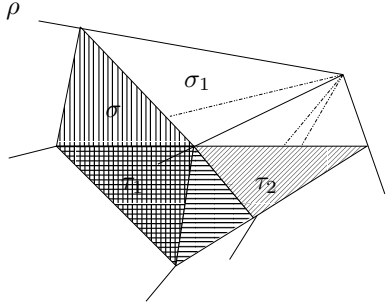
\begin{center}
\input{Fall1.pstex_t}
\end{center}

Note that $\sigma$ and $\tau_1$ are again both $n$--dimensional cones,
whose intersection is also $n$--dimensional.
So in particular, $\sigma\cap \tau_1^{\circ}\ne \emptyset$. 
Moreover, $|\C(\S)^{n}|\geq |\C(S\cap (\sigma+\tau_1))^{n}|$.
In fact, we can reduce the problem to considering 
the pair of cones $\sigma$ and $\tau_1$ in $\S\cap (\sigma+\tau_1)$. 

Because if  $(\sigma+\tau_1)^{\circ}$ is contained in  the support 
of some $f$--equivalence class,
then  it follows from Proposition 
\ref{Kegelabschluss} that the same is true for its
face $\tau':=\sigma_1+\tau_0\subset u^{\perp}$. This implies 
$(\sigma_1+\tau_0)^{\circ}\subset |\sigma_1|_f$, and so by our general
assumption, $\tau'$ must be strictly convex.
It  follows from our choice of $u$ that
$\emptyset\ne u^{\perp}\cap\tau^{\circ}\subset\tau_0^{\circ}\subset(\tau')^{\circ}$,
and therefore $\sigma_1\sim_f \tau$.
Consider the system of fans $\S'$ obtained from $\S\cap(\sigma_1+\tau_2)$
by adding the fan of faces of $\tau'$, and glueing $\tau'$ to the other
cones along the zero cone. Since $\sigma$ lies on the other side of the
hyperplane $u^{\perp}$, we have $|\C(\S')^{n}|<|\C(\S)^{n}|$.
We also know that $\tau_2\cap(\tau')^{\circ}\ne\emptyset$.
It follows by induction that
$(\tau_2+\sigma_1)^{\circ}=(\tau_2+\tau')^{\circ}\subset
|\tau|_f$.
Therefore it suffices to show
the claim for the pair of cones $\sigma$ and $\tau_1$ in
$\S\cap (\sigma+\tau_1)$.

If $\sigma\cup\tau_1$ is not convex, then
we can again decompose $\sigma+\tau_1$ along a defining hyperplane
meeting $\tau_1^{\circ}$ and reduce the problem to considering
the cone $\sigma$ together with a smaller cone $\tau_1'$ as above.
After repeating this procedure a finite number of times however,
we arrive at a pair of cones such that their union is convex.

Now let us consider the second case, namely that $\dim\tau=n-1$ and
$\dim(\tau\cap\sigma)=n-1$, and suppose that
$\tau\not\subset\sigma$. If there is a cone $\sigma'\in\C(\S)^{n}$
with $(\sigma')^{\circ}\cap\sigma^{\circ}\ne\emptyset$ and
$\sigma'\not\subset\sigma$, then we can conclude with the assertion
in the first case applied to the pair of cones $\sigma, \sigma'$
in $\S\cap(\sigma+\sigma')$ that
$(\sigma+\sigma')^{\circ}\subset |\sigma|_f$.
It follows again from the general assumption that
$\sigma'':=\sigma+\sigma'$ is strictly convex.
By replacing the two fans of faces of $\sigma$ and $\sigma'$ respectively
in $\S$ by the single fan of faces of $\sigma''$, where $\sigma''$
is glued to all the other cones only along the zero cone, we obtain
a system of fans $\S'$, that contains strictly fewer $n$--dimensional
cones than $\S$. So by induction the claim follows.

From now on, let us assume that for every cone $\sigma'\in\C(\S)^{n}$
that is not contained in $\sigma$, we have 
$\sigma'\cap\sigma^{\circ}=\emptyset$.
 By assumption, the intersection $\sigma\cap\tau$
contains a point of $\tau^{\circ}$.
Therefore we can find a facet $\sigma_1$ of $\sigma$ with 
$\sigma\cap\tau^{\circ}\ne\emptyset$ and $\sigma_1\not\subset\tau$.
(To see this, note that $\sigma\cap\tau\subsetneq \tau$. 
So there is a facet $\rho_1\prec\sigma\cap\tau$ with 
$\rho_1^{\circ}\subset\tau^{\circ}$.
Since $\sigma=\rho+(\sigma\cap\tau)$, where 
$\rho\not\subset \lin(\sigma\cap\tau)$, $\sigma_1:=\rho+\rho_1$ 
is a facet of $\sigma$ with
 the desired properties.)

As in the previous case, we decompose $\sigma+\tau$ along a defining
hyperplane of $\sigma_1$.
We set 
$$\tau_0:=\tau\cap u^{\perp}, \quad
\tau_1:=\{v\in \tau; u(v)\geq 0\}\quad \text{and}\quad
\tau_2:=\{v\in \tau; u(v)\leq 0\}\,,$$
and observe that $\tau_1$ and $\tau_2$ are cones intersecting
 in the common face
$\tau_0$. And, since $\sigma+\tau=\sigma_1+\tau$,
we again have $\sigma+\tau=(\sigma_1+\tau_1)\cup(\sigma_1+\tau_2)=
(\sigma+\tau_1)\cup (\sigma_1+\tau_2)$.

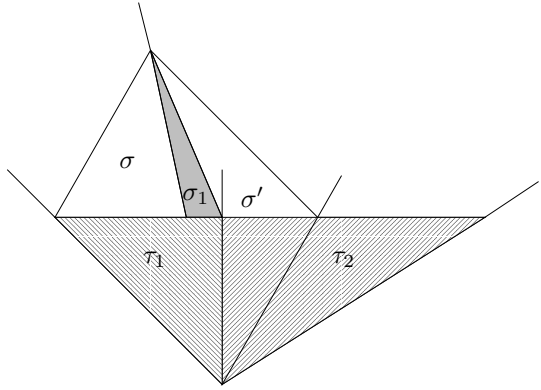
\begin{center}
\input{Fall2.pstex_t}
\end{center}

Moreover, since the $n$--dimensional cones in $\C(\S)$
 cover $\sigma_1$, we can choose a
 cone $\sigma'\in\C(\S)^{n}$ through
a point $v\in\sigma_1\cap\tau^{\circ}$, such that
$\dim(\sigma_1\cap\sigma')=n-1$.
Because of our assumptions on
the position of the $n$--dimensional cones relative to $\sigma$,
that implies in particular that $\sigma'\subset\sigma_1+\tau_2$.

One consequence is that $|\C(\S\cap (\sigma+\tau_1))^{n}|<
|\C(\S)^{n}|$. So by induction we have 
$(\sigma+\tau_1)^{\circ}\subset |\sigma|_f$.
Furthermore, the face $\sigma_1+\tau_0$ of $\sigma+\tau_1$ must
be contained in an $f$--equivalence class. Since 
$\sigma_1\cap \tau_0^{\circ}\ne\emptyset$, that means
$(\sigma_1+\tau_0)^{\circ}\subset |\sigma_1|_f$.
By Lemma \ref{Nachbarkegel},  $\sigma_1\sim_f\sigma$,
and therefore in fact $(\sigma_1+\tau_0)^{\circ}\subset |\sigma|_f$.

Now consider the system of fans $\S'$ obtained from 
$\S\cap(\sigma_1+\sigma')$ by adding the fan of faces of
$\tau':=\sigma_1+\tau_0$ and glueing $\tau'$ to the other cones
along the zero cone.
Since $\sigma'\cap\sigma_1^{\circ}\ne\emptyset$ and
$\dim\sigma_1=\dim\tau$, we have $\sigma'\cap(\tau')^{\circ}\ne\emptyset$.
Because $|\C(\S')^{n}|<
|\C(\S)^{n}|$, we can conclude by induction
that $(\tau'+\sigma')^{\circ}\subset |\sigma'|_f=|\sigma|_f$.
Here the equality again is a consequence of Lemma \ref{Nachbarkegel}.
Altogether we obtain that the relative interior of the convex cone
$$\sigma'':=(\sigma+\tau_1)\cup(\tau'+\sigma')$$
is contained in $|\sigma|_f$.
Replacing the fans of faces of $\sigma$ and $\sigma'$ in $\S$ as
above by the fan of faces of $\sigma''$, glued to the other cones
along zero, we end up with a system of fans $\S''$ that has
strictly fewer $n$--dimensional cones than $\S$. So by induction,
applied to the pair of cones $\sigma''$ and $\tau$,
the claim follows.
This ends the proof.
\endproof

\section{Weakly Proper Quotients}

\noindent
In this section we prove the main result of this paper,
namely that if the toric quotient of a  toric variety with
respect to a subtorus action is weakly proper
 and the quotient variety is of expected dimension, then 
the toric quotient is even categorical.

In order to simplify our arguments, we work in the more
general context of toric prevarieties and
obtain the announced result as a corollary of a statement in this context.

First we note that the construction of toric quotients given
in \cite{Torische Quotienten} also proves the existence of
a separated quotient in the toric category
 for a subtorus action on a toric prevariety. 
More precisely, we have the following:

\begin{proposition} 
Let $X$ be a toric prevariety with big torus $T$, 
and let $H$ be a subtorus of $T$. Then 
there is an $H$--invariant toric morphism $q\colon X\to Y$
to some separated toric variety $Y$ such that
every $H$--invariant toric morphism from $X$ to a toric variety factors
uniquely through $q$.
\end{proposition}

If $X$ is separated then $q$ is precisely the toric quotient of $X$ by $H$. If $X$ is
not separated, we call the morphism $q$ the {\it separated toric quotient\/} of $X$ by $H$.

\bigskip

\proof This statement follows directly from the existence of a quotient
fan of a system of cones proved in \cite{Torische Quotienten}.
Suppose that $X$ arises from an affine system of fans
$\S$  in a lattice $N$, and let $L$ be 
the primitive sublattice of $N$ corresponding to $H$. 
The set $\C(\S)$ of cones of $\S$ satisfies the definition of a system
of $N$--cones given in \cite{Torische Quotienten}. Therefore by
Theorem 2.3 of \cite{Torische Quotienten} there is
a well--defined fan $\Delta$ in a lattice $\t{N}$, the so--called
{\it quotient fan\/} of $\C(\S)$ by $L$, satisfying the following universal property:
Whenever we have a fan $\Delta'$ in a lattice $N'$ and a lattice
homomorphism $F\colon N\to N'$ with $L\subset \ker(F)$, such that
every cone of $\C(\S)$ is mapped into a cone of $\Delta'$,
then there is a unique map of fans from $\Delta$ to $\Delta'$ such
that the diagram commutes. If we translate this property back
into the language of toric prevarieties, then we get back exactly
the universal property stated in the proposition.
\endproof

The main theorem in the context of toric prevarieties is the following:

\begin{theorem}\label{weakly proper toric quotient}
Let $X=X_{\S}$ be a toric prevariety  arising from an affine system of fans $\S$ in
a lattice $N$, and let $H$ be the subtorus corresponding to a given sublattice $L\subset N$.
The quotient fan $\Delta$ of the set $\C(\S)$ of cones of $\S$ is a fan in
a quotient lattice $(N/L)/L'$ of $N/L$. Let $P\colon N\to N/L$ and
$P'\colon N/L\to (N/L)/L'$ denote the projections.

If we assume  $P(|\S|)=(P')^{-1}(|\Delta|)$,
then the separated toric quotient $q\colon X_{\S} \to X_{\Delta}$ with respect to $H$ has
the following universal property (CA): Every $H$--invariant morphism from $X$ to a variety 
factors uniquely through $q$.
\end{theorem}

Let us briefly recall the procedure to
obtain the quotient fan $\Delta$ of $\C(\S)$ by $L$ as described 
in \cite{Torische Quotienten}. 

\medskip

{\it Initialization:\enspace} Let $P\colon N\to N/L$ denote the projection. Let $S_1$ 
denote the set of maximal elements of
$\{ P_{\RR}(\sigma) ; \sigma\in\C(\S)\}$  with respect to  inclusion of sets.

{\it Loop, Step $l>1$:\enspace} If possible choose $\tau,\tau'\in S_{l-1}$ such that 
$\tau\cap\tau'$ is
not a face of $\tau'$ and let $\rho'\prec\tau'$ denote the smallest face 
containing $\tau\cap\tau'$.
Now set $S_{l}$ to be the set of maximal elements of $S_{l-1}\cup \{\tau+\rho'\}$. 
Otherwise set $n:=l-1$
and stop the procedure.

{\it Output:\enspace} Let $\Sigma$ be the set 
of all faces of the cones of $S_n$.
By construction, the set  $\Sigma$ is a {\it quasi-fan\/}, i.e.~the cones in $\Sigma$ may not
be strictly convex, but otherwise all the axioms of a fan are fulfilled. We call $\Sigma$
the quotient quasi-fan of $\C(\S)$ by $L$.

{\it Final step:\enspace} The minimal element $V(\Sigma)$ of  $\Sigma$ 
 is a linear subspace of $(N/L)_{\RR}$, and $V(\Sigma)\cap N/L=L'$ is a sublattice. 
Let $P'\colon (N/L) \to (N/L)/L'$
denote the projection. Then
$\Delta:=\{P'_{\RR}(\tau) ; \tau\in \Sigma\}$.

\bigskip

Our main applications of Theorem~\ref{weakly proper toric quotient} are the following two corollaries.

\begin{corollary}\label{Hauptergebnis}
 For a toric variety $X$ and a subtorus $H$ of the
big torus $T$, let $p\colon X\to Y$ denote
the toric quotient. If $p$ is weakly proper
 and the quotient space $Y$ is of expected dimension, i.e.
$\dim Y=\dim T/H$,
then the toric quotient is a quotient in the category
of algebraic varieties.
\end{corollary}

\proof Suppose that $X$ arises from a fan $\Delta_1$ in a lattice $N$.
Since
$\dim Y= \dim T/H=\dim (N/L)_{\RR}$  the quotient fan 
$\Delta$ of $\Delta_1$ by the lattice $L$ corresponding to $H$
lives in the space $N/L$ and $P'$ is the identity on $N/L$.
Moreover, the projection 
$P\colon N\to N/L$ is the lattice homomorphism associated to the toric quotient $p$.
So if $p$ is weakly proper, then by Proposition~\ref{weakpropchar}, we have
 $P(|\Delta_1|)=|\Delta|$, and Theorem~\ref{weakly proper toric quotient} applies.
\endproof

\begin{corollary}\label{convex case}
Let $X$ be a toric variety corresponding
to a fan with convex support. 
That means that there is a proper toric morphism from 
$X$ onto an affine toric variety. 
Then the toric quotient of $X$
by any subtorus of the big torus is always a quotient in the
category of algebraic varieties.
\end{corollary}

\proof Suppose that $X$ arises from a fan $\Delta_1$ with convex support.
Then any projection from $\Delta_1$ to a quotient quasi-fan is automatically surjective,
since the cones of any quotient quasi-fan are obtained
by successively forming  convex hulls.\endproof

Now we come to the proof of the main theorem.

\bigskip

\noindent{\sc Proof of Theorem 6.2.\enspace}
First we can reduce the situation to a special case.
Let us denote by $\widehat{H}$ the largest subtorus of $T$ such that
every $H$--invariant morphism from $X$ to some variety is in fact 
$\widehat{H}$--invariant. Then in particular, the separated toric quotient
 $q$ is $\widehat{H}$--invariant, and that means that
$H\subset\widehat{H}\subset \ker(q)$. Here as before $\ker(q)\subset T$ denotes
the kernel of the homomorphism of tori induced by $q$.

Note that the separated toric quotients with respect to $H$ and  $\widehat{H}$ 
coincide. Moreover, the assumption $P(|\S|)=(P')^{-1}(|\Delta|)$,
implies that  $\widehat{P}(\vert\S\vert)=(\widehat{P}')^{-1}(|\Delta|)$,
where $\widehat{P}$ denotes the projection modulo the lattice
$\widehat{L}$ corresponding to $\widehat{H}$ and $\widehat{P}'$ the
projection with $P'\circ P= \widehat{P}'\circ \widehat{P}$.
So for this proof we can assume w.l.o.g. that $H=\widehat{H}$.

We claim that moreover we can assume $H=1$ or equivalently $L=0$.
To see this let $f\colon X_{\S}\to Z$ be an $H$-invariant 
morphism to some variety $Z$.
Consider the algebraic quotients of the affine
charts $p_i\colon X_{\sigma_{ii}}\to X_{\sigma_{ii}}\good H$.
Since any good quotient is categorical, the restriction of $f$
to $X_{\sigma_{ii}}$ factors uniquely through $p_i$.
Therefore we obtain morphisms
$\t{f}_i\colon X_{\sigma_{ii}}\good H\to Z$ with $f_i=\t{f}_i\circ p_i$.

In particular, $f$ is $\ker(p_i)$--invariant for every $i$,
and therefore by assumption  $\ker(p_i)=H$ for all $i$.
That means that $X_{\sigma_{ii}}\good H$
is the affine toric variety defined by the cone 
$\tau_i:=P_{\RR}(\sigma_{ii})$ in $N':=N/L$, and
all the cones $\tau_i$ are strictly convex.

If we glue all the affine charts $X_{\tau_i}$ along the open
orbit $T'=T/H$, we obtain a toric prevariety $X'$ that in some
sense is a first non-separated approximation of the quotient variety $X_{\Delta}$.
The prevariety $X'$ corresponds to the system 
of fans  $\S'=(\Delta_{ij}')_{i,j\in I}$ in $N'$, where
$\Delta_{ii}'$ is the fan of faces of $\tau_i$ and $\Delta_{ij}'=\{0\}$
for all $i\ne j$. We can view the quotient fan $\Delta$ of $\C(\S)$ by $L$
as the quotient fan of $\C(\S')$ by the zero lattice. Let
 $q'\colon X_{\S'}\to X_{\Delta}$ denote the separated 
toric quotient of $X'=X_{\S'}$ by the trivial group.

The morphisms $\t{f}_i$ coincide on the open orbit $T'$, and hence
they glue together to a morphism $\t{f}\colon X_{\S'}\to Z$. Since
for every $i$, the morphism $f$ agrees with $\t{f}\circ p_i$ on the affine
chart $X_{\sigma_{ii}}$, we can conclude that
$f$ factors through $q$ if and only if $\t{f}$ factors
through $q'$. 
$$\xymatrix{
X_{\Delta} & X_{\S}  \ar[l]_q \ar[r]^f & Z \cr
& X_{\S'}  \ar[ul]^{q'} \ar[ur]_{\t{f}} & \cr
}
\qquad,$$

Therefore from now on we may assume that $H=\widehat{H}=1$.
In this case the condition in Theorem~\ref{weakly proper toric quotient}
says that $|\S|=|\Sigma|$.

Let $f\colon X_{\S}\to Z$ be a morphism to a variety. 
Let us have a closer look at the algorithmic construction
of the quotient fan $\Delta$ of $\C(S)$ by the zero lattice.
The first set $S_1$ simply consists of the maximal cones in $\C(S)$.
Let $X_1$ denote the toric prevariety obtained from the maximal
affine charts of $X_{\S}$, but now only glued along the open
orbit. The identity on the affine charts induces  a toric morphism 
$p_1\colon X_1\to X_{\S}$. Set $f_1:=f\circ p_1\colon X_1\to Z$.
 
By induction on the  index $l$ counting the steps in the algorithm
we now prove the following for $l>1$:
\begin{enumerate}
\item 
All the cones in $S_l$ are strictly convex.
Hence we may define a toric prevariety $X_l$ from the affine toric varieties
corresponding to the cones in $S_l$ by glueing them along the open orbit $T$.
\item 
There is a morphism 
$f_l\colon X_l\to Z$ such that we have a commutative diagram
$$\xymatrix{
 X_{{l-1}}  \ar[r]^{f_{l-1}} \ar[d]_{p_l} & \ Z \cr
X_{l} \ar[ur]_{f_l} & \cr
}
\qquad,$$
where $p_l$ denotes the toric morphism induced by the identity on $N$.
\end{enumerate}

So let $l>1$, and assume that the induction hypothesis is fulfilled
for all previous steps.
Suppose that $S_l$ was obtained from $S_{l-1}$ by replacing a cone $\tau\in S_{l-1}$ by
$\t\tau=\tau+\rho'$, where $\rho'\prec\tau'\in S_{l-1}$ such that
$(\tau\cap\tau')^{\circ}\subsetneq (\rho')^{\circ}$.
Then in particular, $\tau\cap  (\rho')^{\circ}\ne\emptyset$.

Let $\S_l$ denote the system of fans associated to $X_l$.
By the general assumption, $\vert \S_l\vert=\vert \S_{l-1}\vert$ and hence
 $\t\tau\subset \vert \S_{l-1}\vert$. Since we assumed that $\widehat{H}=1$,
we can apply Proposition \ref{Kegelvergroesserung} to 
$\S_{l-1}\cap \t\tau$ to
 conclude that 
$\t\tau^{\circ}=\tau^{\circ}+(\rho')^{\circ}\subset |\tau|_{f_{l-1}}$. 

So by Proposition \ref{Kegelabschluss}, the relative interior of any
 face of $\t\tau$ is
contained in some $f_{l-1}$--equivalence class. 
If $\t\tau$ would contain a linear subspace, then this would imply
that $f_{l-1}$ for any choice of $f$
would be invariant with respect to the corresponding
subtorus, and that contradicts the assumption that $\widehat{H}=1$.

Therefore in fact $\t\tau$ must be strictly convex, and
$f_{l-1}$ by Corollary \ref{indMor} defines
a morphism on $X_{\t\tau}$. Hence we can extend $f_{l-1}$
to a morphism $f_l\colon X_{\S_l}\to Z$ as desired.
This proves the induction claims (i) and (ii). 

It remains to show that the morphism $f_n\colon X_{\S_n}\to Z$
corresponding to the   last step of the loop, factors through the 
quotient variety $X_{\Delta}$.

As we have just seen, in our special case the cones in $S_n$ are 
strictly convex. Therefore $S_n$ 
is in fact the set of maximal cones of the quotient fan $\Delta$.
The toric prevariety $X_{\S_n}$ is obtained from the  affine charts $X_{\tau}$, where
$\tau\in S_n$, by glueing them along the open orbit $T$.
The morphism $f_n\colon X_{\S_n}\to Z$ induces morphisms
$f_{\tau}\colon X_{\tau}\to Z$ that coincide on the open dense orbit $T$.

We can also view the $X_{\tau}$ as affine charts of $X_{\Delta}$.
But in this variety for any pair $\tau,\tau'\in S_n$,
the charts $X_{\tau}$ and $X_{\tau'}$ are glued along the common
subset $X_{\tau\cap\tau'}$. 
Since $Z$ is separated, and the morphisms $f_{\tau}$ and $f_{\tau'}$
agree on the open dense orbit, they even agree on $X_{\tau\cap\tau'}$.
So they fit together to a morphism $\t{f}\colon X_{\Delta}\to Z$
and $f_n=\t{f}\circ p$, where $p\colon X_{\S_n}\to X_{\Delta}$ is the toric 
morphism induced by the identity on $N$. This ends the proof.
\endproof

We also obtain as a corollary the following statement about actions
of subtori of small codimension
that was proved in \cite{Beispielsammlung} as Theorem 4.1.

\begin{proposition}\label{codim2pre} 
Let  $H\subset T$ be a subtorus of codimension $\leq 2$.
Then for any separated toric quotient  of a toric prevariety
with big torus $T$ with respect to $H$  the universal property~(CA) holds.
\end{proposition}



\noindent{\sc Proof of Proposition \ref{codim2pre}.\enspace} 
Let $p\colon X\to Y$ denote
the separated toric quotient of a toric prevariety $X$ by $H$.
If $\dim T/H=1$, then the conditions of Theorem 
\ref{weakly proper toric quotient}
 are automatically satisfied for $p$,
and therefore the universal property~(CA) holds.

Let us now assume that $\dim T/H=2$.
Then the projections of the maximal
cones of $\Delta$ in $N/L$ are at most $2$--dimensional, and
therefore the convex hull 
of two overlapping projections either equals their
union or equals the whole plane $N_{\RR}/L_{\RR}$.  
So if the    variety $Y$ is $2$--dimensional, 
then the algorithm for constructing 
the quotient fan shows that here $p$ also automatically satisfies the
conditions of the theorem.

If $\dim Y\leq 1$, then among the   cones of the given
toric prevariety $X=X_{\S}$ there must be a chain
of   cones $\sigma_1,\dots,\sigma_r$ such that
$P_{\RR}(\sigma_i^{\circ})\cap P_{\RR}(\sigma_{i+1}^{\circ})\ne\emptyset$,
 $\tau:=\bigcup_{i=1}^{r-1} P_{\RR}(\sigma_i)$ is strictly convex,
and $\tau\cup P_{\RR}(\sigma_r)$ contains a one--dimensional
linear subspace $L'$. Assume that $[\sigma_i,i]\in\Omega(\S)$ and
let $X_i:=X_{[\sigma_i,i]}$ denote the corresponding affine chart in $X$.
Let $X':=\bigcup_{i=1}^{r-1} X_i \subset X$ 
denote the open toric subprevariety of $X$. 
The toric morphism $p'\colon X'\to X_{\tau}$ is weakly proper
and hence (CA) holds for $p'$.

Now let $f\colon X\to Z$ be an $H$--invariant morphism. Then
$f|_{X'}$ factors uniquely through $p'$, i.e.~there is a unique morphism
$\t{f}_1\colon X_{\tau}\to Z$ such that $f=\t{f}_1\circ p'$ on $X'$.
The toric quotient $p''\colon X_{\sigma_r}\to X_{\tau_r}$ 
(where $\tau_r:=P_{\RR}(\sigma_r)$) is even
an affine quotient and hence  categorical. So we similarly obtain 
a unique morphism $\t{f}_2\colon X_{\tau_r}\to Z$ such that
$f=\t{f}_2\circ p''$ on $X_{\sigma_r}$.
The morphisms $\t{f}_1$ and $\t{f}_2$ glue together to a morphism
$\t{f}\colon X_{\S'}\to Z$, where $\S'$ denotes the affine system of fans
obtained from  $\tau$ and $\tau_r$ as maximal cones by glueing them
only along the zero cone.
We have $f=\t{f}\circ q$, where $q\colon X'\cup X_{\sigma_r}\to X_{\S'}$
denotes the natural toric morphism.

Since $\tau^{\circ}\cap \tau_r^{\circ}\ne\emptyset$, we can conclude
from Proposition \ref{Kegelvergroesserung} 
that 
$(\tau + \tau_r)^{\circ}\subset |\tau|_{\t{f}}$ or
there exists a non--trivial linear subspace
$L\subset \tau+\tau_r$ that is contained in the support of an
$\t{f}$--equivalence class. Since $L'\subset (\tau+\tau_r)$,
in any case we have a one--dimensional linear subspace
contained in the support of an $\t{f}$--equivalence class.
and  $\t{f}$ is
invariant with respect to the associated subtorus $H'$ of the
big torus $T'$ of $X_{\S}$. So $f$ is invariant with respect to
$q^{-1}(H')$.
That shows that in fact $p$ is the separated toric quotient of $X$ with
respect to a torus of codimension at most $1$, and we are
back in the first case. 
 \endproof

\section{Example}

In this section we give an example of a $\CC^*$--action on a
$4$--dimensional  toric variety with $3$--dimensional toric
quotient space, such that the toric quotient $p\colon X\to Y$
 is categorical, but not
uniform in the sense of Mumford, i.e. there is a saturated $T$--stable
open subset $U$ of $X$ such that the restriction of $p$ to $U$ is not
the categorical quotient of $U$. Note that \cite{Beispielsammlung}
also contains an example of a categorical quotient  that is not uniform.
But in that example the dimension of the quotient space is not maximal
but strictly less than $\dim(T/H)$.

Let $X=X_{\Delta}$ be the $4$--dimensional toric variety associated
to the fan $\Delta$ in $\ZZ^4$  with the following four
maximal cones:
\begin{eqnarray*}
\sigma_1 &:=&\cone(e_3,e_2+e_3,e_1+e_2+e_3,e_1+e_3) \\
\sigma_2 &:=&\cone(e_1,e_1+e_3,e_1+e_2+e_3,e_1+e_2)  \\
\sigma_3 &:=&\cone(e_2,e_1+e_2,e_1+e_2+e_3,e_2+e_3)  \\
\sigma_4 &:=&\cone(e_1,e_1+e_3-e_2,e_3,e_4) \,.    \\
\end{eqnarray*}
Let $P\colon \ZZ^4\to\ZZ^3$ denote the linear map given by the matrix
$$\left(\matrix{1&0&0&1\cr 0&1&0&0\cr 0&0&1&1\cr
}\right)\,.
$$
Then following diagram shows the images of the maximal cones in
$\ZZ^3$. Note that $P(|\Delta|)$ in this case is a convex cone.

Let $Y$ denote the toric variety associated to the fan in $\ZZ^3$ with
the two maximal cones $\tau_1:=\bigcup_{i=1}^3 P_{\RR}(\sigma_i)$ and
$\tau_2:=P_{\RR}(\sigma_4)$.  The lattice homomorphism $P$ defines a toric morphism $p\colon X\to Y$, and in fact $p$ is the toric quotient of $X$ by
the action of the one--dimensional subtorus $H$ corresponding to 
the kernel of $P$, which is generated by $(1,0,1,-1)^T$ in $(\CC^*)^4$.

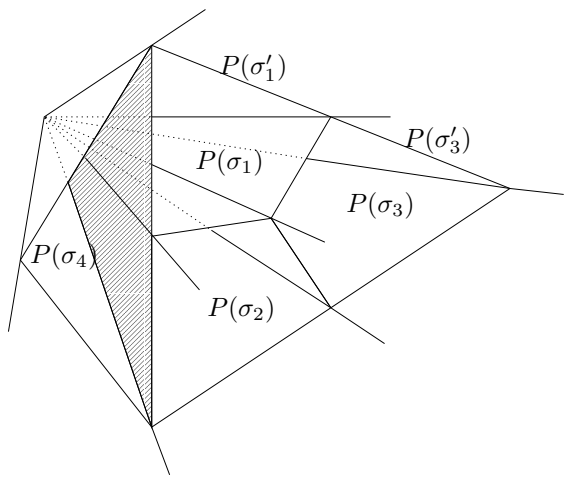
\begin{center}
\input{beispiel.pstex_t}
\end{center}

Since $p$ is weakly proper
and $\dim Y=\dim X-1$, the toric quotient in this case is even
categorical. However the quotient does not satisfy the base--change
property. Consider the open subset 
$U:=X_{\sigma_1'}\cup X_{\sigma_3'}$ of $X$ consisting of the two
affine charts corresponding to the cones $\sigma_1'=\cone(e_3,e_2+e_3)\prec\sigma_1$ and
 $\sigma_3':=\cone(e_2+e_3,e_2)\prec \sigma_3$.
This subset is saturated, and the image under $p$ is the affine toric
subvariety of $Y$ corresponding to the cone generated by
$e_2,e_3$ in $\ZZ^3$.

On the other hand, the toric variety $U$ admits a good quotient, 
and its quotient space 
is $U//H=X_{P(\sigma_1')}\cup X_{P(\sigma_3')}$. So in particular,
$U//H$ is not affine. Therefore  the restriction of $p$ to $U$ 
cannot be the categorical quotient of $U$ by $H$.

\bibliography{}

\begin{thebibliography}{KaStZe}%
%
\small\setlength{\itemsep}{0pt}
\bibitem[AC;Ha, 1]{Torische Quotienten} A.~A'Campo-Neuen, J.~Hausen: Quotients of Toric Varieties by the Action of a Subtorus. T\^ohoku Math. J. {\bf
  51} (1999), 1--12.
%
\bibitem[AC;Ha, 2]{Prevarieties} A.~A'Campo-Neuen, J.~Hausen: Toric Prevarieties
and Subtorus Actions. Geom. Dedicata {\bf 87} (2001), 35-64.
%
\bibitem[AC;Ha, 3]{Beispielsammlung} A.~A'Campo-Neuen, J.~Hausen: 
Examples and Counterexamples for  Existence of Categorical Quotients. 
 J. Algebra, {\bf 231} (2000), 67--85.
%
%
\bibitem[Bo]{Bo} A. Borel: Linear Algebraic Groups, Second Enlarged
  Edition. Sprin\-ger, New York, 1991.
%
%
%
\bibitem[Ew]{Ew} G.~Ewald: Combinatorial Convexity and Algebraic
  Geometry. Sprin\-ger, New York, 1996.
%
\bibitem[Fu]{Fu} W.~Fulton: Introduction to Toric Varieties. Princeton
  University Press, Princeton, 1993.
%
%
%
%
%
\bibitem[Ko]{Ko} J. Koll\`ar: Quotients Spaces Modulo Algebraic Groups.
  Ann.  of Math. {\bf 145} (1997), 33--79.
%
%
\bibitem[Mu]{GIT} D. Mumford: Geometric Invariant Theory. Springer,
  Berlin, 1965.
%
\bibitem[Od]{Od} T.~Oda: Convex Bodies and Algebraic
  Geometry. Springer, Berlin, 1988.
%
\bibitem[Po;Vi]{PV} V.~L.~Popov, E.~B.~Vinberg: Invariant Theory. 
In: Algebraic Geometry IV (A.~N.~Parshin, I.~R.~Shafarevich, eds.), 
Encyclopaedia of Mathematical Sciences {\bf 55}, Springer, Berlin, 1994.
%
%
%
\bibitem[Sw]{Sw} J.~\'Swi\c{e}cicka: Quotients of toric varieties
  by actions of subtori. Colloq. Math., Vol.~82, No.~1 (1999), 105--116.
%
%
\end{thebibliography}

\end{document}